\documentclass[12pt,reqno,tbtags]{amsart}

\date{4 April 2021}

\title[A note on the Markovian SIR epidemic on a random graph with given degrees]
{A note on the Markovian SIR epidemic on a random graph with given degrees}

\author{Malwina Luczak}

\address{School of Mathematics and Statistics, University of Melbourne, 813 Swanston Street, Parkville VIC 3010, Australia.}
\email{mluczak@unimelb.edu.au}


\usepackage{amsmath,amssymb,amsthm,bbm,fullpage} 

\usepackage{ifdraft}
\usepackage{url}
\usepackage[square,numbers]{natbib}
\bibpunct[, ]{[}{]}{;}{n}{,}{,}

\numberwithin{equation}{section}

\renewcommand\le{\leqslant}
\renewcommand\ge{\geqslant}

\allowdisplaybreaks

\theoremstyle{definition}

\theoremstyle{remark}

\newenvironment{romenumerate}[1][0pt]{
\addtolength{\leftmargini}{#1}\begin{enumerate}
 \renewcommand{\labelenumi}{\textup{(\roman{enumi})}}%
 \renewcommand{\theenumi}{\textup{(\roman{enumi})}}%
 }{\end{enumerate}}

\newcounter{oldenumi}
{\setcounter{oldenumi}{\value{enumi}}
\begin{romenumerate} \setcounter{enumi}{\value{oldenumi}}}
{\end{romenumerate}}

\newcounter{thmenumerate}

\newcounter{xenumerate}   

\newcommand\marginal[1]{\marginpar{\raggedright\parindent=0pt\tiny #1}}
\newcommand\REM[1]{{\raggedright\texttt{[#1]}\par\marginal{XXX}}}

\begingroup
  \divide\count255 by 60
  \multiply\count255 by -60
  \advance\count255 by \time
  \ifnum \count255 < 10 \xdef\klockan{\the\count1.0\the\count255}
  \else\xdef\klockan{\the\count1.\the\count255}\fi
\endgroup

\newcommand\set[1]{\ensuremath{\{#1\}}}

\newcommand\lrpar[1]{\left(#1\right)}

\newcommand\abs[1]{|#1|}

\def\rompar(#1){\textup(#1\textup)}    

\newcommand\parfrac[2]{\lrpar{\frac{#1}{#2}}}

\def\xexp(#1){e^{#1}}

\newcommand\downto{\searrow}
\newcommand\upto{\nearrow}
\newcommand\punkt{.\spacefactor=1000}    
\newcommand\iid{i.i.d\punkt}
\newcommand\ie{i.e\punkt}
\newcommand\eg{e.g\punkt}

\newcommand{\tend}{\longrightarrow}

\newcommand\pto{\overset{\mathrm{p}}{\tend}}

\newcounter{CC}
\newcounter{cc}

\renewcommand{\=}{:=}

\usepackage{enumerate}

\numberwithin{equation}{section}

\newcommand{\Prob}{\mathbb{P}}

\newcommand{\N}{\mathbb{N}}
\newcommand{\R}{\mathbb{R}}

\newcommand{\iUp}{T_0} 

\newcommand{\fXS}{h_\mathrm{S}}
\newcommand{\fXI}{h_\mathrm{I}}
\newcommand{\fXR}{h_\mathrm{R}}
\newcommand{\fXX}{h_X}

\newcommand{\fSv}{v_\mathrm{S}}

\newcommand{\pI}{\mathrm p_{\mathrm{I}}}

\newcommand{\aS}{\alpha_{\mathrm{S}}}
\newcommand{\aI}{\alpha_{\mathrm{I}}}
\newcommand{\aR}{\alpha_{\mathrm{R}}}

\newcommand\nS{n_{\mathrm{S}}}
\newcommand\nI{n_{\mathrm{I}}}
\newcommand\nR{n_{\mathrm{R}}}

\newcommand{\nk}{n_{k}}
\newcommand\nSk{n_{\mathrm{S},k}}
\newcommand\nIk{n_{\mathrm{I},k}}
\newcommand\nRk{n_{\mathrm{R},k}}

\newcommand{\mS}{\mu_\mathrm{S}}
\newcommand{\mI}{\mu_\mathrm{I}}
\newcommand{\mR}{\mu_\mathrm{R}}
\newcommand{\mX}{\mu}

\newcommand{\X}[1]{X_{#1}}
\newcommand{\XS}[1]{X_{\mathrm{S},#1}}
\newcommand{\XI}[1]{X_{\mathrm{I},#1}}
\newcommand{\XR}[1]{X_{\mathrm{R},#1}}

\newcommand{\Sv}[1]{S_{#1}}
\newcommand{\Iv}[1]{I_{#1}}
\newcommand{\Rv}[1]{R_{#1}}

\newcommand{\lIv}[1]{\hat I_{#1}}
\newcommand{\lRv}[1]{\hat R_{#1}}

\newcommand{\lTP}[1]{\theta_{#1}}

\newcommand{\Rzero}{\mathfrak{R}_0}

\newcommand{\fXIroot}{\lTP{\infty}}

\newcommand{\calibS}{s_0}

\newcommand\nn{^{(n)}}

\newcommand\Gr{G}
\newcommand\GrCM{G^*}

\begin{document}

\begin{abstract}
This paper provides an additional probabilistic interpretation of the limiting functions for the SIR epidemic on the configuration model derived in the paper by Janson, Luczak and Windridge (2014)~\cite{JLW}.
\end{abstract}

\maketitle

\section{Introduction}
\label{s:intro}

The Markovian SIR process is a simple model for a disease spreading around a finite population
in which each individual is either susceptible, infective or recovered.
Individuals are represented by vertices in a graph $\Gr$ with edges corresponding to potentially infectious
contacts.  Infective vertices become recovered at rate $\rho \ge 0$ and infect each neighbour at rate $\beta > 0$;
those are the only possible transitions, \ie{} recovered vertices never become infective.

The present note concerns SIR epidemics on random graphs with a given degree sequence.
There have been a number of studies of SIR epidemics on random graphs with a given degree sequence~\cite{JLW, JLWH, Volz, Miller, DDMT12, BR}. Here we will give a short discussion and probabilistic interpretation of the results in~\cite{JLW} as well as alternative formulae for the limiting functions derived in that paper.

\section{Model, notation, assumptions and summary of results from~\cite{JLW}}\label{s:notationresults}

Let us start by recalling some notation and assumptions from~\cite{JLW}.

For $n \in \N$ and a sequence $(d_i)_1^n$ of non-negative integers,
let $\Gr = \Gr(n, (d_i)_1^n )$ be a simple graph (\ie{} with no loops or double edges) on $n$ vertices, chosen uniformly at random from among all graphs with degree sequence $(d_i)_1^n$.
(It is tacitly assumed that there is some such graph, so $\sum_{i = 1}^n d_i$ must be even, at least.)

Given the graph $\Gr$, the SIR epidemic evolves as a continuous-time Markov chain.  At any
time, each vertex is either susceptible, infected or recovered.  Each infective vertex recovers at rate $\rho \ge 0$ and also infects each susceptible neighbour at rate $\beta > 0$.

There are initially $\nS$, $\nI$, and $\nR$ susceptible, infective and recovered vertices, respectively.
Further, it is assumed that, for each $k \ge 0$, there are respectively $\nSk$, $\nIk$ and $\nRk$ of these vertices with degree $k$.
Thus, $\nS + \nI + \nR = n$ and $\nS = \sum_{k=0}^\infty \nSk$, $\nI = \sum_{k=0}^\infty \nIk$, $\nR = \sum_{k=0}^\infty \nRk$.
We write $n_k$ to denote the total number of vertices with degree $k$; thus, for each $k$,  $n_k =  \nSk + \nIk + \nRk$.
%

For a sequence $(Y\nn_t)_1^\infty$ of real-valued stochastic processes defined on a subset $E$ of $\R$ and a real-valued function $y$ on $E$,
`$Y\nn_t \pto y(t)$ uniformly on $E' \subseteq E$' means $\sup_{t \in E'}\abs{Y\nn_t - y(t)}
\pto 0$, as $n \to \infty$.

The following regularity conditions for the degree sequence asymptotics are imposed in~\cite{JLW}.

\begin{enumerate}
 \renewcommand{\theenumi}{(D\arabic{enumi})}
 \renewcommand{\labelenumi}{\theenumi}

 \item \label{d:alphaprops} \label{d:first}
The fractions
of initially susceptible, infective and recovered vertices converge to some $\aS, \aI,\aR \in [0,1]$, i.e.
\begin{equation}\label{e:alphaprops}
\nS/n \to \aS, \qquad \nI/n \to \aI, \qquad \nR/n \to \aR.
\end{equation}
Further, $\aS > 0$.

\item \label{d:asympsuscdist}  The degree of a randomly chosen susceptible vertex converges to a probability distribution $(p_k)_{0}^\infty$, i.e.
\begin{equation}\label{e:nSktopk}
\nSk/\nS \to p_k, \qquad k \ge 0.
\end{equation}
Further, this limiting distribution has a finite and positive mean
\begin{equation}\label{e:meansuscdist}
\lambda \= \sum_{k = 0}^\infty k p_k \in (0, \infty).
\end{equation}

\item \label{d:suscmeanUS}
The average degree of a randomly chosen susceptible vertex converges to $\lambda$, i.e.
\begin{equation}
 \label{e:suscmeanUS}
 \sum_{k=0}^{\infty} k \nSk/\nS \to \lambda.
\end{equation}

\item
\label{d:meanIRconv} The average degree over all vertices converges to $\mX > 0$, i.e.
\begin{equation}\label{e:meanX0}
 \sum_{k = 0}^\infty k n_k/n =  \sum_{i = 1}^n d_i/n \to \mX,
\end{equation}
and, in more detail, for some $\mS, \mI, \mR$,
\begin{equation}\label{e:meanSs}
\sum_{k = 0}^\infty k \nSk/n \to \mS,
\end{equation}
\begin{align}\label{e:meanXIR}
  \sum_{k = 0}^\infty k \nIk/n \to \mI,
&\qquad  
  \sum_{k = 0}^\infty k \nRk/n \to \mR.
\end{align}

\item \label{d:SI2ndmoment}
The maximum degree of the initially infective vertices is not too large:
\begin{equation}\label{e:SI2ndmomentx}
\max\set{k:\nIk>0}=o(n).
\end{equation}

\item \label{d:p1orrhopos} Either $p_1 > 0$ or $\rho > 0$ or $\mR > 0$.

\label{d:last}
\end{enumerate}

Clearly, $\aS+\aI+\aR=1$ and $\mS+\mI+\mR=\mX$.  Further, assumptions \ref{d:alphaprops}--\ref{d:suscmeanUS} imply $\sum_{k = 0}^\infty k \nSk/n \to \aS \lambda$, and so
$\mS = \aS \lambda$. 

Let $G^* (n, (d_i)_1^n )$ be the
random multigraph with given degree sequence $(d_i)_1^n$ defined by the
configuration model: we take a set of $d_i$ half-edges for each vertex $i$
and combine half-edges into edges by a uniformly random matching (see \eg{}
\cite{bollobas}). Conditioned on the multigraph being simple, we obtain  $G
= G (n, (d_i)_1^n )$, the uniformly distributed random graph with degree
sequence $(d_i)_1^n$.

Janson, Luczak and Windridge (2014)~\cite{JLW} first prove their results for the SIR epidemic on $\GrCM$, and, by conditioning on $\GrCM$ being simple, they deduce that these results also hold for the SIR epidemic on $\Gr$ .
Their argument relies on the probability that
$\GrCM$ is simple
being bounded away from zero as  $n \to \infty$.
By the main theorem of \cite{Janson:2009} this occurs provided the following condition holds.
\begin{enumerate}
 \renewcommand{\theenumi}{(G\arabic{enumi})}
 \renewcommand{\labelenumi}{\theenumi}
\item\label{d:sumsquaredi=On}
The degree of a randomly chosen vertex has a bounded second moment, i.e.
\begin{equation}\label{e:sumsquaredi=On}
\sum_{k = 0}^\infty k^2 \nk = O(n). 
\end{equation}
\end{enumerate}

The authors of~\cite{JLW} study the SIR epidemic on the multigraph $G^*$, revealing its edges dynamically while the epidemic spreads.
The process analysed in~\cite{JLW} works as follows.
A half-edge is said to be {\it free} if it is not yet paired to another half-edge. A half-edge is called susceptible,
infective or recovered according to the type of vertex it belongs to.

At time $0$, there are $d_i$ half-edges attached to vertex $i$, for each $i$, and all half-edges are free.
Subsequently, each free
infective half-edge chooses a free half-edge at rate $\beta$, uniformly at random from among all the other free half-edges.
Together the pair form an edge, and are removed from the pool of free
half-edges.
If the chosen half-edge belongs to a susceptible vertex then that vertex becomes infective, and thus all of its half-edges become infective also.
Infective vertices also recover at rate $\rho$.

The process stops when there are no free
infective half-edges,  at which point
the epidemic stops spreading.  Some infective vertices may remain but
they will recover at \iid{}  exponential times without affecting
any other vertex, and are irrelevant from of the point of view of the epidemic.
Some susceptible and recovered half-edges may also remain, and these are
paired off uniformly at time $\infty$
to reveal the remaining edges in $\GrCM$.
This step is unimportant for the spread of the epidemic, but is performed for the purpose of transferring the results from the multigraph $\GrCM$ to the simple graph $G$.

Clearly, if all the pairings are completed then the resulting graph is the multigraph $G^*$.
Moreover, the quantities of interest (numbers of susceptible, infective and recovered vertices at each time $t$) have the same distribution as if we were to reveal the multigraph $G^*$ first and run the SIR epidemic on $G^*$ afterwards.

For $t \ge 0$, let $\Sv{t}$, $\Iv{t}$ and $\Rv{t}$ denote the numbers of susceptible,
infective and recovered vertices, respectively, at time $t$.
Thus $\Sv{t}$ is decreasing and $\Rv{t}$ is increasing.
Also $\Sv{0} = \nS$, $\Iv{0} = \nI$ and $\Rv{0} = \nR$.

For the dynamics described above (with half-edges paired off dynamically, as described), for $t \ge 0$, let $\XS{t}$, $\XI{t}$ and $\XR{t}$ be the number of free susceptible, infective and recovered half-edges
at time $t$, respectively.  Thus $\XS{t}$ is decreasing,
$\XS{0} = \sum_{k = 0}^\infty k \nSk$,  $\XI{0} = \sum_{k = 0}^\infty k
\nIk$ and $\XR{0} = \sum_{k = 0}^\infty k\nRk$.
%

For the uniformly random graph $G$ with degree sequence $(d_i)_1^n$, the variables $\XS{t}$, $\XI{t}$ and $\XR{t}$, for $t \ge 0$, are defined as above conditioned on the final multigraph $G^*$ being a simple graph.

It is shown in~\cite{JLW} that, upon suitable scaling, the processes $\Sv{t}, I_t, R_t$, $\XS{t}, \XI{t}, \XR{t}$
converge to deterministic functions.  The limiting functions
are written in terms of a parameterisation
$\theta_t \in [0,1]$ of time solving an ordinary differential equation
given below.  In~\cite{JLW},
the function $\theta_t$ is interpreted as the limiting probability that a given
initially susceptible half-edge has been selected for pairing with a (necessarily infective) half-edge
by time $t$.
Let
\begin{equation}\label{e:fSv}
 \fSv(\theta) \= \aS \sum_{k=0}^\infty p_k \theta^k, \qquad \theta \in [0,1],
\end{equation}
so the limiting fraction of susceptible vertices is $\fSv(\theta_t)$ at time $t$ (since the events of being selected for pairing will be approximately independent for different half-edges, when $n$ is large).  Similarly, for susceptible half-edges, the limiting function is
\begin{equation}\label{e:fXS}
  \fXS(\theta) \=  \aS \sum_{k=0}^\infty k \theta^k p_k = \theta\fSv'(\theta), \qquad \theta \in [0,1].
\end{equation}
For the total number of free half-edges, let
\begin{equation}
\fXX(\theta) := 
\mX \theta^2, \qquad \theta \in [0,1].
\label{e:fX}
\end{equation}
For the numbers of half-edges of the remaining types, for $\theta \in [0,1]$, let
\begin{align}
  \fXR(\theta) &:= \mR \theta + \frac{\mX\rho}{\beta}\theta(1-\theta), \label{e:fXR} \\
  \fXI(\theta) &:= \fXX(\theta) - \fXS(\theta) - \fXR(\theta). \label{e:fXI}
\end{align}
Thus $\fXX(\theta) = \fXS(\theta)+\fXI(\theta)+ \fXR(\theta)$.
Note that
\begin{align}
&&\fSv(1)&=\aS, &
\\ \label{e:fX1}
\fXS(1)&=\aS\lambda=\mS,&
\fXR(1)&=\mR,&
\fXI(1)&=\mX-\mS-\mR=\mI.
\end{align}

The following is shown in~\cite{JLW}.

When $\mI > 0$,  there is a unique $\fXIroot \in (0,1)$ with $\fXI(\fXIroot) = 0$.  Further,
$\fXI$ is strictly positive on $(\fXIroot,1]$
 and strictly negative on $(0,\fXIroot)$.
Defining the `infective pressure'
\begin{equation}\label{e:pI}
\pI(\theta) \= \frac{\fXI(\theta)}{\fXX(\theta)},
\end{equation}
there is a unique solution
$\lTP{t}:[0,\infty) \to (\fXIroot,1]$ to the differential equation
\begin{equation}\label{e:dlTPtmIg0}
\frac{d}{dt} \lTP{t} = -\beta \lTP{t} \pI(\lTP{t}),
\end{equation}
subject to the initial condition $\lTP{0} = 1$.

Furthermore, there is a unique solution $\lIv{t}$  to
\begin{equation}\label{e:dlIvtTmIg0}
 \frac{d}{dt} \lIv{t} = \frac{\beta \fXI(\lTP{t})\fXS(\lTP{t}) }{\fXX(\lTP{t})} - \rho \lIv{t}, \; t \ge 0,\qquad \lIv{0} = \aI.
\end{equation}

Defining also $\lRv{t} \= 1 - \fSv(\lTP{t}) - \lIv{t}$,
Theorem 2.6 in~\cite{JLW} states that, for the epidemic on the multigraph $\GrCM$, under conditions \ref{d:first}--\ref{d:last},
uniformly on $[0,\infty)$,
\begin{align}
\label{e:convmIg0SIR}
\Sv{t}/n &\pto \fSv(\lTP{t}), &
 \Iv{t}/n &\pto \lIv{t}, & \Rv{t}/n &\pto \lRv{t},
\\
\label{e:convmIg0X}
\XS{t}/n &\pto \fXS(\lTP{t}), & \XI{t}/n &\pto \fXI(\lTP{t}),
&\XR{t}/n &\pto \fXR(\lTP{t}),
\\
\label{e:total}
\X{t}/n & \pto \fXX(\lTP{t}).
\end{align}
Moreover, the number $\Sv{\infty} \= \lim_{t \to \infty} \Sv{t}$ of susceptibles that escape infection satisfies
\[
\Sv{\infty}/n \pto \fSv(\fXIroot).
\]

The same holds on the graph $\Gr$ under the additional assumption \ref{d:sumsquaredi=On} (Theorem 2.7 in~\cite{JLW}).


\bigskip

For the case where there are initially a small number of infectives (of order less than $n$, so that $\mI = 0$),
we recall from~\cite{JLW}
\begin{equation}\label{e:R0}
\Rzero \= \parfrac{\beta}{\rho + \beta} \parfrac{\aS}{\mX}
  \sum_{k=0}^\infty (k-1) k p_k;
\end{equation}
the {\it basic reproductive ratio} of the epidemic.
It is shown in~\cite{JLW} that, when $\Rzero > 1$, even if $\mI = 0$, then
there is a unique $\fXIroot \in (0,1)$ with
$\fXI(\fXIroot) = 0$, and that $\fXI$ is strictly positive on $(\fXIroot,1)$
 and strictly negative on $(0,\fXIroot)$.

The initial condition of the limiting differential equation,
now defined on $(-\infty, \infty)$, is shifted
so that $t = 0$ corresponds to the time $\iUp$ in the random process, which is the infimum of times $t$ such that
the fraction of susceptible individuals has fallen from about $\aS= \fSv (1)$ to some
fixed smaller $\calibS$ by time $t$. It is shown in~\cite{JLW} that there is a unique continuously differentiable $\lTP{t}:\R \to (\fXIroot,1)$
such that
\begin{equation}\label{e:dlTPt}
\frac{d}{dt} \lTP{t} = -\beta \lTP{t} \pI(\lTP{t}),
\qquad \lTP{0} = \fSv^{-1}(\calibS).
\end{equation}
Furthermore, $\lTP{t}\downto \fXIroot$ as $t\to\infty$ and $\lTP{t}\upto 1$ as $t\to -\infty$.

The processes are extended to be defined on
$(-\infty,\infty)$ by taking $\Sv{t} = \Sv{0}$ for $t < 0$, and similarly
for the other processes.

The following is proven in~\cite{JLW} (Theorems 2.9 and 2.10), for both the simple graph $\Gr$ and the multigraph $\GrCM$.
Suppose that conditions \ref{d:first}--\ref{d:last} and \ref{d:sumsquaredi=On} hold. Assume that $\Rzero > 1$. Suppose also that $\aI=\mI = 0$ but there is initially at least one infective vertex with non-zero degree.

Then, $\liminf_{n \to  \infty} \Prob(\iUp < \infty) > 0$.
Also, conditional on $\iUp < \infty$,
uniformly on $(-\infty,\infty)$,
\begin{align}
\Sv{\iUp+t}/n &\pto \fSv(\lTP{t}), &
\Iv{\iUp+t}/n &\pto \lIv{t}, & \Rv{\iUp+t}/n &\pto \lRv{t},
\label{e:convscSIR} \\
\XS{\iUp+t}/n &\pto \fXS(\lTP{ t}),& \XI{\iUp+t}/n &\pto \fXI(\lTP{t}), &
\XR{\iUp+t}/n &\pto \fXR(\lTP{t}),
\label{e:convscX}
\\
\label{e:total-1}
\X{\iUp+t}/n & \pto \fXX(\lTP{t}).
\end{align}

Also, conditional on $\iUp < \infty$, the number of susceptibles
that escape infection
satisfies
\[
\Sv{\infty}/n \pto \fSv(\fXIroot).
\]

Here, again, $\lIv{t}$ is the unique solution to
\begin{equation}\label{e:dlIvtT}
 \frac{d}{dt} \lIv{t} = \frac{\beta \fXI(\lTP{t})\fXS(\lTP{t}) }{\fXX(\lTP{t})} - \rho \lIv{t}, \qquad \lim_{t \to -\infty} \lIv{t} = 0,\\
\end{equation}
and $\lRv{t} \= 1 - \fSv(\lTP{t}) - \lIv{t}$.

\section{New probabilistic interpretation of $\theta$ and alternative formulae for limiting deterministic functions}
\label{s:theta-interpretation}

We will now give a more complete probabilistic interpretation of the function $\theta_t$ used to define the deterministic limit for the SIR epidemic.

As stated in the introduction, the function $\theta_t$ used to define the deterministic limits
satisfies
$$\frac{d \theta_t}{dt} = - \beta \theta_t \frac{h_I (\theta_t)}{h_X (\theta_t)}.$$
Substituting $\fXX (\theta) = \mu \theta^2$, $\fXI (\theta) = \mu \theta^2 - \mu_R \theta - \frac{\mu \rho}{\beta} \theta (1-\theta) - \aS \sum_k k p_k \theta^k$, we can rewrite this as
\begin{eqnarray*}
\frac{d \theta_t}{dt} & = & - \beta  \frac{\mu \theta_t^2 - \mu_R \theta_t - \frac{\mu \rho}{\beta} \theta_t (1-\theta_t) - \alpha_S \sum_k k p_k \theta_t^k}{\mu \theta_t}\\
& = & - (\beta + \rho ) \theta_t + \rho + \frac{\beta \mu_R}{\mu} + \frac{\beta \alpha_S}{\mu} \sum_k k p_k \theta_t^{k-1}.
\end{eqnarray*}
It follows that
\begin{eqnarray*}
\frac{d}{dt} (\theta_t e^{(\beta + \rho) t}) = \Big ( \rho + \frac{\beta \mu_R}{\mu} \Big ) e^{(\beta + \rho ) t} + \frac{\beta \alpha_S}{\mu} (\sum_k k p_k \theta_t^{k-1}) e^{(\beta + \rho ) t},
\end{eqnarray*}
and so, integrating,
\begin{eqnarray*}
\theta_t & = & \Big ( 1 - \frac{ \rho + \frac{\beta \mu_R}{\mu}}{\beta + \rho} \Big ) e^{-(\beta + \rho ) t} +  \frac{ \rho + \frac{\beta \mu_R}{\mu}}{\beta + \rho} + \frac{\beta \alpha_S}{\mu} \sum_k k p_k \int_0^t \theta_s^{k-1} e^{-(\beta + \rho )(t-s)} ds\\
& = & \frac{\rho}{\beta + \rho} + \frac{\beta}{\beta + \rho} \frac{\mu_R}{\mu} + \frac{\beta}{\beta + \rho} \frac{\mu - \mu_R}{\mu} e^{-(\beta + \rho) t} \\
& + & \frac{\beta \alpha_S}{\mu} \sum_k k p_k \int_0^t \theta_s^{k-1} e^{-(\beta + \rho )(t-s)} ds,
\end{eqnarray*}
and so
\begin{eqnarray}
\label{eq-formula-theta-0}
\theta_t & = & \frac{\mu_R}{\mu} + \frac{\mu - \mu_R}{\mu} F(t) + \frac{\beta \alpha_S}{\mu} \sum_k k p_k \int_0^t \theta_s^{k-1} e^{-(\beta + \rho )(t-s)} ds,
\end{eqnarray}
where
$$F(t) =  \frac{\rho}{\beta + \rho} + \frac{\beta}{\beta + \rho} e^{-(\beta + \rho) t}.$$

Noting that
\begin{eqnarray*}
\int_0^t \theta_s^{k-1} \frac{dF(t-s)}{ds} ds + \int_0^t \frac{d}{ds} (\theta_s)^{k-1} F(t-s) ds = \Big [\theta_s^{k-1} F(t-s) \Big ]_0^t,
\end{eqnarray*}
we see that, for each $k \ge 1$,
\begin{eqnarray*}
\beta \int_0^t \theta_s^{k-1} e^{-(\beta + \rho) (t-s)} ds & = & \theta_t^{k-1} - F(t)
-\int_0^t \frac{d}{ds} (\theta_s)^{k-1} F(t-s) ds\\   
& = & \theta_t^{k-1} - F(t) + \beta (k-1) \int_0^t \theta_s^{k-1} \frac{\fXI(\theta_s)}{\fXX (\theta_s)} F(t-s) ds.
\end{eqnarray*}
This then implies, using $\aS \sum_k k p_k = \mS$, that
\begin{eqnarray*}
\theta_t & = & \frac{\mu_R}{\mu} + \frac{\mu - \mu_R}{\mu} F(t) + \frac{\alpha_S}{\mu} \sum_k k p_k \theta_t^{k-1} - \frac{\mS}{\mu} F(t)\\
& + & \frac{\beta \alpha_S}{\mu} \sum_k k (k-1) p_k \int_0^t \theta_s^{k-1} \frac{\fXI(\theta_s)}{\fXX (\theta_s)} F(t-s)  ds \\
\end{eqnarray*}
and hence that
\begin{eqnarray}
\label{eq-formula-theta}
\theta_t & = & \frac{\mu_R}{\mu} + \frac{\mI}{\mu} F(t) + \frac{\alpha_S}{\mu} \sum_k k p_k \theta_t^{k-1} \\
& + & \frac{\beta \alpha_S}{\mu} \sum_k k (k-1) p_k \int_0^t \theta_s^{k-1} \frac{\fXI(\theta_s)}{\fXX (\theta_s)} F(t-s)  ds. \nonumber
\end{eqnarray}

Considering formula~\eqref{eq-formula-theta}, we will now discuss how the function $\theta_t$ is the asymptotic probability that a half-edge does not transmit infection (i.e. initiate a pairing) by time $t$. This should be the same as the limiting probability that a given
initially susceptible half-edge has not been paired with a (necessarily infective) half-edge
by time $t$, as interpreted in~\cite{JLW}, since that probability is that its eventual partner has not transmitted infection by time $t$.

Given a random half-edge, conditional on it being initially recovered, which has probability $\mu_R/\mu$, it does not transmit infection by time $t$ with probability 1.

Conditional on the half-edge being initially infected, which has probability $\mu_I/\mu$, it does not transmit by time $t$ with probability $F(t)$. In the formula for $F(t)$, the term $\frac{\rho}{\beta + \rho}$ is the probability that recovery of the vertex occurs before the half-edge's clock goes off. The term $\frac{\beta}{\beta + \rho} e^{-(\beta + \rho) t}$ is the probability that the clock of the half-edge goes off before recovery but neither of these events happens by time $t$.

Conditional on the half-edge being initially susceptible, which happens with probability $\mu_S/\mu$, we need to further consider the degree of its vertex. With conditional probability $\frac{\aS k p_k}{\mS}$, it has degree $k$, and then the edge cannot transmit if the vertex does not get infected by time $t$ or only gets infected by transmitting the infection to the half-edge itself, which happens with probability $\theta_t^{k-1}$. The half-edge also cannot transmit by time $t$ if one of the other $k-1$ half-edges gets infected at some time $s \le t$, but then the clock of the half-edge in question does not go off before vertex recovery over a period of length $t-s$; this happens with probability $-\int_0^t \frac{d}{ds} (\theta_s)^{k-1} F(t-s) ds = \beta (k-1) \int_0^t \theta_s^{k-1} \frac{\fXI (\theta_s)}{\fXX (\theta_s)} F(t-s) ds$.

Alternatively, we have
\begin{eqnarray*}
n \mu \theta_t^2
& = & n \mu_R \theta_t  + n \mI \theta_t F(t) + n \fXS (\theta_t)\\
& + & n \beta \alpha_S \theta_t \sum_k k (k-1) p_k \int_0^t \theta_s^{k-1} \frac{\fXI(\theta_s)}{\fXX (\theta_s)} F(t-s)  ds.
\end{eqnarray*}

The left hand-side here is approximately the total number of free half-edges at time $t$. The term $n \mu_R \theta_t$ is approximately the total number of initially recovered half-edges that are still free at time $t$. The term $n \mI \theta_t F(t)$ is approximately the total number of initially infective half-edges that are still free at time $t$. The term $n h_S (\theta_t)$ is approximately the total number of free susceptible half-edges at time $t$. The term
$$n \beta \aS \theta_t \int_0^t \sum_k k (k-1) p_k \theta_s^{k-1} \frac{\fXI (\theta_s)}{\fXX (\theta_s)} F(t-s) ds$$
is approximately the total number of half-edges belonging to initially susceptible vertices that got infected before time $t$ and are still free at time $t$.

The function $\theta_t$ is closely related to the corresponding function in~\cite{Miller,Volz}, but these papers do not engage with the construction of the configuration model multigraph and simple graph by pairing half-edges.

\medskip

We saw in Section~\ref{s:intro} (equation~(\ref{e:fXR})) that $\XR{t}/n$ is asymptotically close to
\begin{eqnarray*}
 \fXR(\theta_t) &:= \mR \theta_t + \frac{\mX\rho}{\beta}\theta_t(1-\theta_t). 
\end{eqnarray*}
The term $\frac{\mX\rho}{\beta}\theta_t(1-\theta_t)$ in the above formula is compact but does not appear readily interpretable.

We claim that the limiting function can instead be expressed in the form
\begin{eqnarray}
\label{eq-recovered-half-edges}
\tilde{h}_R (t) & = & \mu_R \theta_t + \mu_I \theta_t \frac{\rho}{\beta + \rho} ( 1 - e^{-(\beta + \rho) t} )\\
& + & \aS \sum_k k p_k \theta_t  \Big ( - \int_0^t \frac{d}{ds} (\theta_s)^{k-1} \frac{\rho}{\beta + \rho} (1-e^{-(\beta + \rho )(t-s)}) ds \Big ).
\nonumber
\end{eqnarray}
To understand this formula, note that $n \mR \theta_t$ is approximately the number of free recovered half-edges that were initially recovered.

Also,
$$\frac{\rho}{\beta + \rho} ( 1 - e^{-(\beta + \rho) t} )$$
is the probability that a vertex infectious at time $0$ recovers by time $t$ and that its recovery happens before the clock of an infectious half-edge attached to this vertex goes off. This implies that
$$n \mu_I \theta_t \frac{\rho}{\beta + \rho} ( 1 - e^{-(\beta + \rho) t} )$$
is approximately the total number of free recovered half-edges that were infectious at time $0$.

Finally,
$$n \aS \sum_k k p_k \theta_t  \Big ( - \int_0^t \frac{d}{ds} (\theta_s)^{k-1} \frac{\rho}{\beta + \rho} (1-e^{-(\beta + \rho )(t-s)}) ds \Big )$$
is approximately the total number of free recovered half-edges whose vertices were susceptible at time $0$, got infected and recovered by time $t$.

We are now going to verify that $\tilde{h}_R(t) = \fXR (\theta_t)$.
This means that we need to verify that
\begin{eqnarray*}
\frac{\mu \rho}{\beta} (1-\theta_t) & = & \mu_I \frac{\rho}{\beta + \rho} ( 1 - e^{-(\beta + \rho) t} )\\
& + & \alpha_S \sum_k k p_k   \Big ( - \int_0^t \frac{d}{ds} (\theta_s)^{k-1} \frac{\rho}{\beta + \rho} (1-e^{-(\beta + \rho )(t-s)}) ds \Big ).
\end{eqnarray*}

To do that, first note that, integrating by parts,
\begin{eqnarray*}
& - \int_0^t \frac{d}{ds} (\theta_s)^{k-1} \frac{\rho}{\beta + \rho} (1-e^{-(\beta + \rho )(t-s)}) ds\\
& = \Big [ - \theta_s^{k-1} \frac{\rho}{\beta + \rho} (1-e^{-(\beta + \rho) (t-s)}) \Big ]^t_0
- \rho \int_0^t \theta_s^{k-1}  e^{-(\beta + \rho) (t-s)} ds\\
& =  \frac{\rho}{\beta + \rho} (1-e^{-(\beta + \rho) t}) -   \rho \int_0^t \theta_s^{k-1}  e^{-(\beta + \rho) (t-s)} ds.
\end{eqnarray*}

This means we actually need to verify that
\begin{eqnarray*}
\frac{\mu \rho}{\beta} (1-\theta_t) & = & \mu_I \frac{\rho}{\beta + \rho} ( 1 - e^{-(\beta + \rho) t} )\\
& + & \mu_S \frac{\rho}{\beta + \rho} (1-e^{-(\beta + \rho) t}) \\
& - & \rho \alpha_S \sum_k k p_k \int_0^t \theta_s^{k-1} e^{-(\beta + \rho) (t-s)} ds.
\end{eqnarray*}

But, as seen in~\eqref{eq-formula-theta-0},
\begin{eqnarray*}
\theta_t
& = & \frac{\mu_R}{\mu} + \frac{\mu - \mu_R}{\mu} F(t) + \frac{\beta \alpha_S}{\mu} \sum_k k p_k \int_0^t \theta_s^{k-1} e^{-(\beta + \rho )(t-s)} ds,
\end{eqnarray*}
where
\begin{eqnarray*}
F(t) = \frac{\rho}{\beta + \rho} + \frac{\beta}{\beta + \rho} e^{-(\beta + \rho) t},
\end{eqnarray*}
and so

\begin{eqnarray*}
- \rho \alpha_S \sum_k k p_k \int_0^t \theta_s^{k-1} e^{-(\beta + \rho) (t-s)} ds & = & -\frac{\mu \rho}{\beta} \theta_t  + \frac{\rho (\mu   - \mu_R)}{\beta + \rho}  e^{-(\beta + \rho ) t}\\
& + &  \frac{\rho}{\beta} \frac{ \rho \mu }{\beta + \rho} + \frac{\rho \mu_R}{\beta + \rho}.
\end{eqnarray*}
This means that we need to verify that
\begin{eqnarray*}
\frac{\mu \rho}{\beta} (1-\theta_t) & = & \mu_I \frac{\rho}{\beta + \rho} ( 1 - e^{-(\beta + \rho) t} ) + \mu_S \frac{\rho}{\beta + \rho} (1-e^{-(\beta + \rho) t})\\
& - & \frac{\mu \rho}{\beta} \theta_t + \frac{\rho (\mu   - \mu_R)}{\beta + \rho}  e^{-(\beta + \rho ) t} + \frac{\rho}{\beta} \frac{ \rho \mu }{\beta + \rho} + \frac{\rho \mu_R}{\beta + \rho},
\end{eqnarray*}
which holds, noting that $\mu_S + \mu_I + \mu_R = \mu$.

\bigskip

Similarly, we have an alternative formula for the limit of $\XI{t}/n$, the asymptotic scaled number of free infectious half-edges at time $t$:
\begin{eqnarray}
\label{eq-formula-infectious-half-edges}
\tilde{h}_I (t) & = &  \mu_I \theta_t e^{-(\beta + \rho) t}
 +  \aS \theta_t \sum_k k p_k   \Big ( - \int_0^t \frac{d}{ds} (\theta_s)^{k-1}  e^{-(\beta + \rho )(t-s)} ds \Big ).
\end{eqnarray}

For infectious vertices, we have
\begin{eqnarray}
\label{eq-formula-infectious-vertices}
\lIv{t} & = & \aI e^{-\rho t} + \beta \int_0^t \fXI (\theta_s) \frac{\fXS (\theta_s)}{\fXX (\theta_s)} e^{-\rho (t-s)} ds \\
& = &  \aI e^{-\rho t} - \int_0^t \frac{d  \fSv(\theta_s)}{ds} e^{-\rho (t-s)} ds, \nonumber
\end{eqnarray}
and, for recovered vertices,
\begin{eqnarray}
\label{eq-formula-recovered-vertices}
\lRv{t} & = & \aR + \aI (1 - e^{-\rho t}) + \beta \int_0^t \fXI (\theta_s) \frac{\fXS (\theta_s)}{\fXX (\theta_s)} (1-e^{-\rho (t-s)}) ds \\
& = &  \aR + \aI (1 - e^{-\rho t}) - \int_0^t \frac{d  \fSv(\theta_s)}{ds} (1-e^{-\rho (t-s)} ds). \nonumber
\end{eqnarray}

\end{document}